\documentclass[reqno,11pt]{amsart}
\usepackage{amsmath, latexsym, amsfonts, amssymb, amsthm, amscd}
\usepackage{graphics,epsf,psfrag}
\numberwithin{equation}{section}

\newtheorem{theorem}{Theorem}[section]
\newtheorem{lemma}[theorem]{Lemma}

\newtheorem{rem}[theorem]{Remark}

\DeclareMathOperator{\sign}{\mathrm{sign}}

\newcommand{\ind}{\mathbf{1}}

\newcommand{\R}{\mathbb{R}}
\newcommand{\Z}{\mathbb{Z}}
\newcommand{\N}{\mathbb{N}}
\renewcommand{\tilde}{\widetilde}

\newcommand{\cH}{{\ensuremath{\mathcal H}} }

\newcommand{\cN}{{\ensuremath{\mathcal N}} }
\newcommand{\cL}{{\ensuremath{\mathcal L}} }

\newcommand{\cD}{{\ensuremath{\mathcal D}} }

\newcommand{\bP}{{\ensuremath{\mathbf P}} }
\newcommand{\bE}{{\ensuremath{\mathbf E}} }

\DeclareMathSymbol{\leqslant}{\mathalpha}{AMSa}{"36} % nicer `smaller or equal'
\DeclareMathSymbol{\geqslant}{\mathalpha}{AMSa}{"3E} % nicer `larger or equal'
\DeclareMathSymbol{\eset}{\mathalpha}{AMSb}{"3F}     % nicer `emptyset'
\newcommand{\dd}{\,\text{\rm d}}             % a straight d for differentials
\newcommand{\maxtwo}[2]{\max_{\substack{#1 \\ #2}}} % max with 2 lines
\newcommand{\mintwo}[2]{\min_{\substack{#1 \\ #2}}} % min with 2 lines
\newcommand{\sumtwo}[2]{\sum_{\substack{#1 \\ #2}}} % sum with 2 lines
\newcommand{\prodtwo}[2]{\prod_{\substack{#1 \\ #2}}}     % product 2 lines
\newcommand{\bbE}{{\ensuremath{\mathbb E}} }

\newcommand{\bbL}{{\ensuremath{\mathbb L}} }

\newcommand{\bbP}{{\ensuremath{\mathbb P}} }

\newcommand{\ga}{\alpha}
\newcommand{\gb}{\beta}
\newcommand{\gga}{\gamma}            % \gg already exists...
\newcommand{\gd}{\delta}
\newcommand{\go}{\omega}
\newcommand{\gto}{{\tilde\omega}}

\newcommand{\gS}{\Sigma}

\makeatletter
\def\captionfont@{\footnotesize}
\def\captionheadfont@{\scshape}

\long\def\@makecaption#1#2{%
  \vspace{2mm}
  \setbox\@tempboxa\vbox{\color@setgroup
    \advance\hsize-6pc\noindent
    \captionfont@\captionheadfont@#1\@xp\@ifnotempty\@xp
        {\@cdr#2\@nil}{.\captionfont@\upshape\enspace#2}%
    \unskip\kern-6pc\par
    \global\setbox\@ne\lastbox\color@endgroup}%
  \ifhbox\@ne % the normal case
    \setbox\@ne\hbox{\unhbox\@ne\unskip\unskip\unpenalty\unkern}%
  \fi
  \ifdim\wd\@tempboxa=\z@ % this means caption will fit on one line
    \setbox\@ne\hbox to\columnwidth{\hss\kern-6pc\box\@ne\hss}%
  \else % tempboxa contained more than one line
    \setbox\@ne\vbox{\unvbox\@tempboxa\parskip\z@skip
        \noindent\unhbox\@ne\advance\hsize-6pc\par}%
\fi
  \ifnum\@tempcnta<64 % if the float IS a figure...
    \addvspace\abovecaptionskip
    \moveright 3pc\box\@ne
  \else % if the float IS NOT a figure...
    \moveright 3pc\box\@ne
    \nobreak
    \vskip\belowcaptionskip
  \fi
\relax
}
\makeatother
\def\writefig#1 #2 #3 {\rlap{\kern #1 truecm
\raise #2 truecm \hbox{#3}}}

\newcommand{\tf}{\textsc{f}}

\newcommand{\rf}{\mathtt f}

\newcommand{\cone}{\mathtt{C}}

\begin{document}
\title[Smoothing of the depinning transition]{
Smoothing effect of quenched disorder\\ 
on polymer depinning  transitions}

\author{Giambattista Giacomin}
\address{Laboratoire de Probabilit{\'e}s de P 6\ \& 7 (CNRS U.M.R. 7599)
  and  Universit{\'e} Paris 7 -- Denis Diderot,
U.F.R.                Mathematiques, Case 7012,
                2 place Jussieu, 75251 Paris cedex 05, France
\hfill\break
\phantom{br.}{\it Home page:}
{\tt http://www.proba.jussieu.fr/pageperso/giacomin/GBpage.html}}
\email{giacomin\@@math.jussieu.fr}
\author{Fabio Lucio Toninelli}
\address{
Laboratoire de Physique, ENS Lyon (CNRS U.M.R. 5672), 46 All\'ee d'Italie, 
69364 Lyon cedex 07, France
\hfill\break
\phantom{br.}{\it Home page:}
{\tt http://perso.ens-lyon.fr/fabio-lucio.toninelli}}
\email{fltonine@ens-lyon.fr}
\date{\today}

\begin{abstract}
We consider general disordered
models of 
pinning of directed  polymers on a defect line. 
This class contains in particular the  $(1+1)$--dimensional 
interface wetting model, the disordered Poland--Scheraga model of DNA denaturation
and other $(1+d)$--dimensional polymers in interaction with flat interfaces.
We consider also the case  of copolymers with adsorption at a selective interface.
  Under quite 
general conditions, these models are known
to have a (de)localization transition at some critical line
in the phase diagram.  In this work we prove in particular  that,
as soon as disorder is present, the transition is at least of second order, in the sense
that the free energy is differentiable at the 
critical line, so that the order parameter vanishes continuously at the transition. 
On the other
hand, it is known that the corresponding non--disordered models can have
a first order (de)localization transition, with a discontinuous
first derivative. Our result shows therefore 
that the presence of the disorder has really
a smoothing effect on the transition.
The relation with the predictions based on the Harris criterion is discussed.
\end{abstract}

\maketitle

\section{Introduction and models}

Quenched disorder is  expected  to smooth 
phase transitions in many situations. This is for instance the
case of ferromagnetic spin systems in dimension $d\le 2$ (if the spins 
are discrete)
or in dimension $d\le 4$ (if  they have rotation symmetry):  
when a  random magnetic field is present the  Imry--Ma
argument \cite{cf:Imry-Ma}, made rigorous by Aizenman and Wehr
\cite{AizWeh90}, implies that
these models do not 
exhibit, at any temperature, the first order phase transition with
associated spontaneous magnetization  which characterizes the corresponding
pure models. For the analogous phenomenon 
for SOS effective interface models 
 see \cite{cf:BK}.

In the present work, under some conditions on the disorder distribution, we prove that a similar effect takes place 
in  models of directed polymers in random media
exhibiting a localization/delocalization transition on a defect line.
Such a transition may be of first or higher order in the corresponding
pure, i.e. non--disordered, cases and we show that,  as soon as 
disorder is present, 
the transition is at least of second order.
It is  important to emphasize that the mechanism inducing
the smoothing of the transition in our case is very different from
the Imry--Ma one, and it is rather based on an estimate of the 
probability that the polymer visits rare but very favorable regions 
where the disorder produces a large positive fluctuation of the
partition function. 

We will consider mainly two classes of models: random pinning (or wetting) models
and random copolymers at selective interfaces.  
In pinning models \cite{cf:AS,cf:Derrida,cf:Fogacs,cf:Nicolas,cf:TangChate} 
the typical situation one has in mind is that of a directed path in $(1+d)$ dimensions, which receives 
a reward (or a penalty) at each intersection with a $1$--dimensional region (a {\sl defect line}, in the 
physical language) according to whether 
the charge present on the line at the  intersection point is positive or negative. 
On the other hand the copolymer model, whose study was initiated in \cite{cf:GHLO} in the theoretical physical literature
(see \cite{cf:Monthus} and references therein for updated physics developments) and in 
 \cite{cf:BdH,cf:Sinai} in the mathematical one, aims
at modeling a $(1+1)$--dimensional, directed heteropolymer containing both hydrophobic and hydrophilic components,
in presence of an interface (the line $S=0$ in the notations of Section \ref{sec:model_c}) 
separating two solvents, situated in the upper and lower half--planes ($S>0$ and $S<0$).
One solvent favors the hydrophilic components and the other 
favors the hydrophobic ones, i.e., if the charge of the $n$--th monomer is positive (negative), 
this monomer tends to 
be in the upper (lower) half--plane. 

The main question one would like to answer is whether or not the interaction  induces a localization 
of the polymer along the defect line or interface.

\subsection{Random pinning and wetting models}
\label{sec:model_p}
Let $S:=\left\{ S_n \right\}_{n=0,1, \ldots}$ be a homogeneous
process
on an arbitrary set that contains a point $0$, with $S_0=0$ and law $\bP$.
For $\gb\ge0, h\in \mathbb R$ and $\go=\{\go_n\}_{n=1,2,\cdots}\in \R^\N$,
one introduces the probability measure
\begin{equation}
\label{eq:RNS}
\frac{\dd \bP _{N, \go}}{\dd \bP} (S) \, =\,\frac1{Z_{N,\go}}
\exp\left( \sum_{n=1}^N \left(\gb \go_n -h\right)\ind _{S_n=0}\right) \ind _{ S_N=0}.
\end{equation} 
The choice of
setting $S_N=0$ is just for technical convenience, see Remark \ref{rem:free} below.
Let us set $\tau_0=0$ and, for $i\in \N:= \{1,2, \ldots\}$,  $\tau_{i}= \inf\left\{n> \tau_{i-1}:\, S_n=0 \right\}$ if
 $\tau_{i-1} < +\infty$
and $\tau_i =+\infty$ if  $\tau_{i-1} = +\infty$. We assume that
 $\tau:=\left\{ \tau_i \right\}_i$ is the sequence of partial sums 
of an IID sequence of random variables taking values 
in $\N  \cup \{ \infty\}$ with discrete density $K(\cdot)$:
\begin{equation}
\label{eq:ret_times}
K(n)=\bP(\tau_1=n).
\end{equation}
This is of course the case if $S$ is a Markov chain and for definiteness one should keep this case in mind.
In order to avoid trivialities, we assume that $K(\infty ) <1$ and for the model to be defined 
we assume also that there exists $s\in \N$ such that 
$\bP (S_N=0)>0$ for every $N\in  s\N$. Therefore, starting
with \eqref{eq:RNS},  we always implicitly 
assume that  $N\in  s\N$, when $N$ is the length of the polymer.

Moreover the cases we will consider are such that
there exist $\ga \ge 1$ such that
\begin{eqnarray}
\label{eq:RW}
\liminf_{ n \to \infty} \frac{\log K(sn)}{\log n} \ge -\ga, 
\end{eqnarray}
 and $K(n)=0$ if  $n\not \in s\N$.
 We advise the reader who feels uneasy with the weak requirement \eqref{eq:RW}
 to focus on the case of $K(n)$ behaving like $n^{-\ga}$, possibly times a {\sl slowly varying function}
 (see Appendix~\ref{app:b0}).
 Note that, for $\alpha< 2$, 
the return times to zero  of $S$ are not integrable.
We stress that we have introduced $s$ to account for the possible periodicity of  
 $S$: of course the most natural example is that
 of the simple random walk on $\Z$, $\bP(S_i-S_{i-1}=1)=\bP(S_i-S_{i-1}=-1)=
1/2$ and $\{S_i-S_{i-1}\}_i$ IID, for which  \eqref{eq:RW} holds with $s=2$ and 
$\alpha=3/2$ \cite[Ch. III]{cf:Feller1}. Another example is 
that of the simple random walk on $\Z^d$, $d\ge 2$: in this case, 
$\alpha=d/2$.

We can rewrite $Z_{N,\go}$, and in fact the model itself,
in terms of the sequence $\tau$: the partition function in \eqref{eq:RNS} is
\begin{equation}
Z_{N, \go} \, :=\, \bE\left[ \exp \left( \cH_{N, \go}(S)\right); \, \tau_{\cN _N}=N\right],
\label{eq:withtau}
\end{equation}
with $\cN_N := \sup\{i: \tau_i \le N\}$ and 
\begin{equation}
\label{eq:Hpinning}
 \cH_{N, \go}(S)\, =\,  \sum_{i=1}^{\cN_N}
\left(\gb \go_{\tau_i} -h\right).
\end{equation}

The  disordered pinning model, which we consider here, 
is obtained  assuming that the sequence $\go$ is chosen as a typical realization of an IID
sequence of random variables with law $\bbP$, still denoted by $\go = \{ \go_n \}_n$. 
We assume  finiteness of exponential moments: 
\begin{eqnarray}
  \bbE[\exp(t\go_1)]< \infty,
\end{eqnarray}
for $t\in\R$  and, without loss of generality, $\bbE \,\go_1=0$ and  $\bbE \,\go_1^2=1$.
Further assumptions on $\bbP$ will be formulated in Section \ref{sec:smooth}, where we state our main results.

Under the above assumptions on the disorder
the {\sl quenched free energy} of the model exists, namely
the limit
\begin{equation}
\label{eq:free_energy_p}
\tf(\gb,h)\, := \, \lim_{N\to \infty} \frac 1N \log  Z_{N,\go},
\end{equation}
exists  $\bbP (\dd \go)$--almost
surely and in the $\bbL^1\left( \bbP\right)$ sense.  
The existence of this limit
can be proven via standard super--additivity arguments based on Kingman's sub--additive
ergodic theorem \cite{cf:King} (we refer for example to \cite{cf:AS,cf:BdH,cf:G} for details).
This approach yields also automatically the so called {\sl self--averaging} property
of the free energy, that is the fact that $\tf(\gb,h)$ is non random.
 A simple but fundamental observation is that
\begin{equation}
\label{eq:delocfe_p}
\tf(\gb,h) \ge \, 0. 
\end{equation} 
The proof of such a result is elementary:
\begin{eqnarray}
\label{eq:step_deloc}
\frac 1N \bbE \log  Z_{N,\go}  &\ge& 
\frac 1N \bbE \log \bE 
\left[ 
\exp\left( \sum_{i=1}^{\cN_N}
\left(\gb \go_{\tau_i} -h\right)
\right)
;  \tau_1=N
\right]
\\
&=&\,
\frac 1N { \left(\gb\bbE\left[\go_N \right]-h\right) } 
+
\frac 1N 
\log \bP \left( \tau_1=N\right)\, \stackrel{N \to \infty}{\longrightarrow}\, 0,
\end{eqnarray}
where 
we have  used 
assumption \eqref{eq:RW}.
The proof of \eqref{eq:delocfe_p} suggests the following
partition of the parameter space (or {\sl phase diagram}):
\smallskip
\begin{itemize}
\item The localized region: $\cL = \left\{ (\beta,h) : \, \tf(\beta,h)> 0\right\}$;
\item The delocalized region:  $\cD = \left\{ (\beta,h): \, \tf(\beta,h)=0\right\}$.
\end{itemize}
\smallskip

We set $h_c(\gb) := \sup\{h: \tf (\gb, h)>0\}$, and we will
call $h_c(\gb)$ the {\sl critical point}.
Since $\tf (\gb, \cdot)$ is not increasing and continuous, $\cD =\{(\gb, h):\, h\ge h_c(\gb)\}$.
It is rather easy to see that $|h_c(\gb)|<\infty$, see e.g. \cite{cf:AS} and \cite{cf:G}.

\medskip

\begin{rem}\rm
\label{rem:free}
It would be of course as natural to consider the model with partition function
$Z_{N, \go}^{\rf}:=\bE\left[ \exp \left( \cH_{N, \go}(S)\right)\right]$.
 It is therefore worth to stress that
 by standard arguments, see e.g.  \cite{cf:G}, one sees that 
  \begin{equation}
\begin{split}
  Z_{N, \go} \, \le\,  Z_{N, \go}^{\rf} \, &\le\,  \maxtwo{n \in s\N:}{n \le N} \frac{\sum_{j=n}^\infty K(j)}{K(n)}
 \exp\left(\left\vert \gb \go_N -h \right\vert \right) Z_{N, \go} \\
& \,\le\,
\frac1{\mintwo{n \in s\N:}{n \le N} K(n)}  \exp\left(\left\vert \gb \go_N -h \right\vert \right) Z_{N, \go}
 ,
\end{split}
 \end{equation}
 uniformly in $N \in s \N$ and $\go$. In particular, by \eqref{eq:RW},  the free energy is unaffected by the
 presence of 
 the constraint $\tau_{\cN_N}=N$: we stick to the constrained case because 
it simplifies some technical steps of the proofs.
\end{rem}
\medskip

\begin{rem}\rm
For ease of exposition we have used the generic denomination of {\sl pinning model},
but our framework, 
as it is possibly  clearer  when the partition function
is cast in the form \eqref{eq:withtau}, 
 includes a variety of models, like for example the Poland--Scheraga
model of DNA denaturation, for which theoretical arguments on 
models with excluded volume interactions suggest a value
of $\ga$ larger than $2$ \cite{cf:KMP}:
note that whether the transition in the disordered case 
is of first or higher order is
a crucial and controversial issue in the field, see for example
 \cite{cf:CH,cf:C,cf:GM2}. 
We stress also that, since we can choose $K(\infty)>0$, it is rather easy 
to see that also the disordered $(1+1)$--dimensional  interface wetting
models \cite{cf:Derrida,cf:Fogacs} enter the general class
we are considering.  
\end{rem}

It is  known that the nonrandom case $\gb=0$ is exactly solvable  and that,
according to the law of $\tau_1$, the transition at $h_c(0)=\log(1-K(\infty))$ 
can be either of first or higher order. For completeness and to match our set--up,
we give a quick self--contained analysis of this case in Appendix \ref{app:b0}.

\subsection{Random copolymers at a selective interface}
\label{sec:model_c}
In the case of the copolymer model, the natural setting is to assume, in addition, that the 
state space of the process $\{S_n\}_n$ is $\Z$, that the law $\bP$ is invariant under the 
transformation $S\to -S$ and that $(S_{i}-S_{i-1})\in\{-1,0,+1\}$ for every
$i\ge 1$. 
For instance, if $\{S_{i}-S_{i-1}\}_i$ is a sequence of IID
variables then it is a classical result~\cite[Ch. XII.7]{cf:Feller2} that  $K(n)= c(1+o(1))n^{-3/2}$, for large values of $n \in s\N$ (so $\ga =3/2$): 
$c$ is a constant that can be expressed
in terms of $\bP(S_1-S_0=0)$ and
of
course $s=1$ unless $\bP(S_1-S_0=0)=0$. 
The copolymer model is defined introducing
\begin{equation}
\label{eq:Hcop}
\cH_{N, \go}(S)\, =\, \frac12
\sum_{n=1}^N(\beta \go_n+h)\sign (S_n),
\end{equation}
with
 the convention that $\sign(0)=+1$, and the corresponding partition function
\begin{equation}
Z_{N, \go} \, :=\, \bE\left[ \exp \left( \cH_{N,\go}(S)\right); \, S_N=0\right].
\end{equation}
We may assume without loss of generality that both $\gb $ and $h$ are nonnegative.
The factor $1/2$ in \eqref{eq:Hcop} is introduced just for convenience 
(see Section \ref{sec:prova_cop}).
As in Section \ref{sec:model_p}, the corresponding random model is obtained by choosing $\go$   
as a realization of 
an IID sequence of centered random variables of unit variance and finite exponential moments. 
In analogy with  \eqref{eq:free_energy_p} and \eqref{eq:delocfe_p}, the limit
\begin{equation}
\label{eq:free_energy_c}
f(\gb,h)\, := \, \lim_{N\to \infty} \frac 1N \log  Z_{N,\go},
\end{equation}
exists  $\bbP (\dd \go)$--almost
surely and in the $\bbL^1\left( \bbP\right)$ sense and one can prove as in \eqref{eq:step_deloc} that
$\tf(\gb,h):=f(\gb,h)-h/2\ge 0$.
Therefore, also in this case we can partition the phase diagram into a localized and a delocalized phase,
as
\smallskip
\begin{itemize}
\item The localized region: $\cL = \left\{ (\beta,h) : \, \tf(\beta,h)> 0\right\}$;
\item The delocalized region:  $\cD = \left\{ (\beta,h) : \, \tf(\beta,h)= 0\right\}$.
\end{itemize}
Also in this case, 
we set $h_c(\gb) := \sup\{h: \tf (\gb, h)>0\}$ and we observe that
$\cD =\{(\gb,h) : h\ge h_c(\gb)\}$.

\medskip
Many results have been proven for copolymers at selective interfaces, 
but they are almost always about the case of simple random walks.
However they can be extended in a rather straightforward way to the general case we consider here
 (we omit the details also because
they are not directly pertinent to the content of this paper).
Above all we have that 
 $0< h_c (\gb)< \infty$ for every $\gb >0$. Even more, the explicit bounds
 carry over to the general $S$ we consider here: the upper bound   
  \cite{cf:BdH} is even independent of the choice of the law of $S$, while the lower bound
  \cite{cf:BG} depends on $\ga$.
  
  \medskip
We take this opportunity to stress also that
various results about path behavior in  the two regions are available, both for the 
pinning problem and for the copolymer, and, once again,
mostly in the simple random walk case
 $\alpha=3/2$. For instance, it is known that in $\cL$ long excursions of the walk from the line
$S=0$ are exponentially suppressed and  the fraction of sites where the polymer crosses the line remains nonzero in the 
thermodynamic limit
(e.g., cf \cite{cf:G,cf:Sinai} and references therein, and \cite{cf:GTloc}). On the other hand,
 in the interior of $\cD$ the number of intersections is, in a suitable sense,  $O( \log N)$ \cite{cf:GT} and, for the copolymer, the number of steps where $\sign(S_n)=-1$
is also   $O( \log N)$ \cite{cf:GT}. 

\section{Smoothing of the depinning transition}
\label{sec:smooth}
While the free energy can be proven to be infinitely differentiable with respect to all of its parameters
 in the
region $\cL$ \cite{cf:GTloc}, no results are available about its regularity at the critical point,
apart from the obvious fact that $\tf$ is continuous since it is convex.
The result of the present paper 
partly fills this gap, showing that the transition is at least of second order, as soon as disorder is present.

In order to state our main theorem, we need some further assumptions on the disorder variables $\omega$.
We will consider two distinct cases:
\medskip

\newcounter{Lcount}
\begin{list}{{\bf C\arabic{Lcount}}:}
{\usecounter{Lcount}
 \setlength{\rightmargin}{\leftmargin}}
\item {\it Bounded random variables}. The random variable $\go_1$ is bounded,
\begin{equation}
  \label{eq:omegabound}
  |\go_1|\le M<\infty.  
  \end{equation}
  \item {\it Unbounded continuous random variables}. 
The law of $\omega_1$ has a  density $P(\cdot)$ with respect to the Lebesgue measure on $\R$, and there exists 
$0<R<\infty$ such that
\begin{eqnarray}
\label{eq:exentropia}
\int_{\R} P(y+x)\log \left(\frac{P(y+x)}{P(y)}\right) \dd y\le R x^2
\end{eqnarray}
in a neighborhood of  $x=0$. 
This is true in great generality whenever $P(\cdot)$ is positive,
for example when the disorder is Gaussian and, more generally, whenever 
 $P(\cdot)=\exp(-V(\cdot))$, with $V(\cdot)$  a polynomial bounded below.
\end{list}
\medskip

Then, one has:
\label{C1C2}
\begin{theorem}
\label{th:smooth_cont} 
Under condition {\bf C1} or {\bf C2}, both for the copolymer and for the pinning
model, for every $0<\beta<\infty$  
there exists $0< c(\beta)<\infty$, possibly depending on $\bbP$, 
such that for every $1\le \alpha<\infty$ 
  \begin{eqnarray}
    \label{eq:smooth_cont}
    \tf(\beta,h)\le \alpha c(\beta) (h_c(\beta)-h)^2
  \end{eqnarray}
if $h<h_c(\beta)$.
\end{theorem}

Although the above result, coupled of course with \eqref{eq:delocfe_p}, seems to be in the same spirit as the rounding effect
proven by Aizenman and Wehr \cite{AizWeh90} for the two--dimensional Random Field Ising Model, the physical mechanisms of 
smoothing are deeply different in the two cases.
While \cite{AizWeh90}  is based on a rigorous version of the Imry-Ma argument \cite{cf:Imry-Ma}
(i.e., a comparison between the effect of boundary conditions and of disorder fluctuations in the bulk
due to the random magnetic field)
in our case the boundary conditions play no role at all and everything is based on 
an energy--entropy 
argument inspired by \cite{cf:BG}.

\begin{rem}\rm
  It is important to observe that, as explained in Appendix \ref{app:b0},  in the 
    pinning case 
the deterministic model, $\beta=0$,
has a {\sl first order phase transition} whenever  $\sum_{n\in \N}n K(n)<+\infty$
 \cite{cf:AS}, in particular when $\ga >2$. Theorem \ref{th:smooth_cont}  
therefore shows that the disorder has really a smoothing effect on the transition.
But more than that is true: for $\ga \in (3/2,2)$, 
$\tf (0,h)$ at $h_c(0)$ is strictly less regular than $\tf(\gb, h)$ at $h_c(\gb)$.
In fact $\tf (0,h) > (h_c(0)- h)^{2-\gd}$ for  $\gd \in (0, (2\ga -3)/(\ga-1))$ and small values of $h_c(0)-h>0$ 
(for sharper results, see Appendix \ref{app:b0}).
Notice that this is in agreement with the so--called {\sl Harris criterion} \cite{cf:Harris}, which
predicts that arbitrarily weak disorder modifies the nature of a second--order phase transition as soon as the 
critical exponent of the specific heat in the pure case is positive. In the present situation, this 
condition corresponds just to $\alpha>3/2$. The Harris criterion also predicts that the critical behavior does not
change if $\alpha<3/2$, which is compatible with Theorem \ref{th:smooth_cont}.
Rigorous work connected to the Harris criterion, in the Ising model context,
may be found in \cite{cf:CCFS}.

As one realizes easily from the proof of Theorem \ref{th:smooth_cont} given in Section \ref{sez:prova}, 
the  constant $c(\beta)$  in 
\eqref{eq:smooth_cont} can be very large (of order $O(\beta^{-2})$) for $\beta$ small. This is rather
intuitive: for $\beta\to0$ one approaches the deterministic situation, where the transition can be of first order.

\end{rem}

\begin{rem}\rm
In the theoretical physics literature, the 
(de)localization transition is claimed to be in some cases of order higher than two \cite{cf:TM}, or even
of infinite order \cite{cf:TangChate, cf:Monthus}. The method we present here, which is rather insensitive to the details of the model,
does not allow to prove more than second order in general. It is likely that  finer results 
require model--specific
techniques. 
\end{rem}

\begin{rem}\rm
\textit{A generalized model: copolymers with adsorption}.
It is also possible to consider copolymer models with an additional pinning interaction \cite{cf:SW}, 
as
\begin{equation}
\label{eq:Hcop+p}
\cH_{N, \go,\gto}(S)\, =\,  \sum_{n=1}^{N}
\left(\gb_1 \go_n -h_1\right)\ind_{S_n=0}+\frac12
\sum_{n=1}^N(\beta_2 \tilde \go_n-h_2)\sign (S_n).
\end{equation}
Here, both $\go$ and $\gto$ are sequences of IID centered  random variables with unit variance and 
finite exponential moments. In addition, one assumes $\go$ to be independent from $\gto$.
This model corresponds to the situation where the interface between the two solvents is not neutral.
While for simplicity we will not present details for this model, we  sketch here what happens in this case.
In analogy with the previous models, one can partition the phase diagram (i.e., the space of the parameters 
$\gb_{1},h_{1},\gb_2,h_2$) 
into a localized and a delocalized region, separated by a {\sl critical surface}. 
In this case, Theorem \ref{th:smooth_cont} is easily generalized to 
give that the free energy $\tf(\gb_1,h_1,\gb_2,h_2)$ has continuous first derivatives with respect to $h_{1},h_2$ 
when these parameters approach the critical surface from the localized region.
\end{rem}

\section{Proof of the smoothing effect}
\label{sez:prova}

\subsection{The pinning case}
In this section we prove Theorem \ref{th:smooth_cont} for the 
pinning case, and in the next one we explain how the proof can be immediately extended to the 
copolymer model.

The key idea, in analogy with \cite{cf:GT}, is to introduce a new free energy where the fraction of 
sites where the polymer comes back to zero is fixed. In other words, recalling that
\begin{equation}
  \label{eq:enne}
  \cN_{N}=|\{1\le n\le N: S_n=0\}|
\end{equation}
one introduces, for $m\in[0,1]$,
\begin{equation}
  \label{eq:def_phi}
  \begin{split}
  \phi(\beta,m)\, &=\, 
  \lim_{\varepsilon\searrow 0}\lim_{N\to\infty}\frac1N \bbE \log \hat Z_{N,\go}(\beta;m,\varepsilon)\\
&:=\, \lim_{\varepsilon\searrow 0}\lim_{N\to\infty}\frac1N \bbE \log \bE \left(e^{\beta\sum_{n=1}^N\go_n\ind_{S_n=0}}
\ind_{(\cN_N/N)\in[m-\varepsilon,m+\varepsilon]}\ind_{S_N=0}\right).
\end{split}
\end{equation}
Note that the limit is well defined, since $\bbE \log \hat Z_{N,\go}(\beta;m,\varepsilon)$ is super--additive in 
$N$, thanks to the IID assumption on the increments of
the sequence of return times $\tau$, and non--increasing for  $\varepsilon\searrow 0$. 
Therefore the limit $\phi_\varepsilon (\gb, m)$ of $(1/N)\bbE \log \hat Z_{N,\go}(\beta;m,\varepsilon)$, 
as $N \to \infty$, exists, as well as the second limit $\phi_\varepsilon (\gb, m)\searrow \phi (\gb, m)$
as $\varepsilon \searrow 0$. Notice moreover that \eqref{eq:def_phi}
holds also without taking the expectation, as for $(1/N)\log Z_{N,\go}$: in this case of course 
the limit $N\to \infty$ has to be taken in the $\bbP (\dd \go)$--a.s. sense.
Moreover, it is immediate to realize that
$\phi(\beta,0)=0$ and that 
\begin{eqnarray}
\label{eq:ULB}
 \phi(\beta,m)\le \tf(\beta,0)\le \beta,
\end{eqnarray}
so that $\phi(\beta,m)$ is always bounded above. Finally, always thanks 
the IID property of the differences of successive return times to zero, it is easy to show
that $\phi(\beta,\cdot)$ is concave:
\begin{eqnarray}
  \label{eq:phi_conc}
  \phi(\beta,m)\ge x \phi(\beta,m_1)+(1-x)\phi(\beta,m_2)
\end{eqnarray}
if $0\le x\le 1$ and $m=x m_1+(1-x)m_2$.
By exploiting the  $\bbP(\dd\go)$--a.s. convergence of $(1/N) \log \hat Z_{N,\go}(\beta;m,\varepsilon)$
to the nonrandom limit $\phi_\varepsilon (\gb,m)$ and the subsequent convergence
for $\varepsilon \searrow 0$, one deduces
that $\tf(\beta,h)$ and $\phi(\beta,m)$ are related by a Legendre transform:
\begin{eqnarray}
  \label{eq:Legendre}
  \tf(\beta,h)=\sup_{m\in[0,1]}\left(\phi(\beta,m)-h m\right).
\end{eqnarray}
In turn, this allows to identify $h_c(\beta)$ in terms of $\phi(\beta,\cdot)$ as
\begin{eqnarray}
  \label{eq:mu_c}
  h_c(\beta)=\inf\left\{h : \sup_{m\in[0,1]}\left(\phi(\beta,m)-h m\right)=0\right\}.
\end{eqnarray}
The key technical step in the proof of Theorem \ref{th:smooth_cont} is the following:
\begin{theorem}
\label{th:stimaphi}
Under condition {\bf C1} or {\bf C2}, for every $0<\beta<\infty$ 
there exists  $0<\cone(\beta)<\infty$ such that, for every $1\le\alpha<\infty$,
  \begin{eqnarray}
    \label{eq:stimaphi}
    \phi(\beta,m)- h_c(\beta)m\le -\frac{\cone(\beta)}\alpha m^2 
  \end{eqnarray}
if $0\le m\le 1$.
\end{theorem}
\medskip

\noindent {\it Proof of Theorem  \ref{th:smooth_cont}}. It is an immediate consequence of 
Theorem \ref{th:stimaphi} and  \eqref{eq:Legendre}.
In fact, 
by  \eqref{eq:stimaphi} we have
\begin{eqnarray}
\label{eq:locale}
    \phi(\beta,m)- h m\le -\frac{\cone(\beta)}\alpha m^2+(h_c(\gb)-h)m, 
  \end{eqnarray}
for every $m\in[0,1]$.  
Taking the supremum over $m$ on both sides of 
 \eqref{eq:locale}, by  
 \eqref{eq:Legendre} we obtain 
 \begin{eqnarray}
   \label{eq:21}
\tf(\gb,h)&\le& \sup_{m\in [0,1]}
\left(-\frac{\cone(\beta)}\alpha m^2+(h_c(\gb)-h)m\right)\\\nonumber
&\le& \sup_{m\ge 0}\left(-\frac{\cone(\beta)}\alpha m^2+(h_c(\gb)-h)m\right)=\frac\alpha{4\cone(\beta)}(h_c(\gb)-h)^2.
 \end{eqnarray}

\hfill $\stackrel{\text{Theorem  \ref{th:smooth_cont}, pinning}}{\Box}$

\medskip

We go now to the proof of Theorem \ref{th:stimaphi}. We will consider first
the case in which $\go_1$ satisfies condition {\bf C2}, because it 
technically lighter. The two cases differ only in the first part of the proof
(that is, up to Remark~\ref{rem:Hineq} below),
where the probability of a rare event is estimated from below 
by changing the law of the disorder and by evaluating the corresponding
relative entropy price. In the case of  {\bf C2} it is sufficient to {\em shift} 
(i.e. to {\em translate}) the distribution of the 
disorder variables, and the relative entropy estimate implied by \eqref{eq:exentropia}
fits well the rest of the proof.
Under assumption {\bf C1}, instead, we have to {\em tilt} the law of $\go$
and the arising expressions need to be re--worked, see Lemma \ref{lemma_psi} below, 
before stepping to the second part of the proof. 

%In the case of  {\bf C2} the change of measure is obtained 
%by adding a suitable constant
%to the disorder variables, we call such an operation {\sl shift}, and the estimate is
%automatically cast in a form that fits the rest of the proof.
%Under assumption {\bf C1} instead we  use the standard {\sl tilting} procedure
%and the arising expressions need to be re--worked, see Lemma \ref{lemma_psi} below, 
%before stepping to the second part of the proof. 

\medskip

\noindent
{\it Proof of Theorem \ref{th:stimaphi} under assumption {\bf C2}}.
Due to the concavity of $\phi(\beta,\cdot)$, it is enough to prove \eqref{eq:stimaphi} for $0< m\le c_1$, with
$c_1=c_1( \beta)>0$, not depending on $\alpha$.
We define
\begin{eqnarray}
  \label{eq:assurdo}
  \gga(\beta,m):=-\phi(\beta,m)+h_c(\beta)m+c_2\frac{\beta^2 m^2}{\alpha}>0,
\end{eqnarray}
 where $c_2>0$ is a constant depending only 
on $\bbP$, which will be chosen later.  
We stress that the term containing $c_2$ has been added
simply because a priori one knows only that $-\phi(\beta,m)+h_c(\beta)m\ge 0$,
cf. \eqref{eq:mu_c}, and it turns out to be technically practical to
work with $\gga(\beta,m)>0$.
For $\ell \in s\N$ we define also
\begin{eqnarray}
  \label{eq:evento}
  A_{\ell,m,\varepsilon}=
\left\{\omega\in \R^\ell: \log \hat Z_{\ell,\go}(\beta;m,\varepsilon)- h_c(\beta)m\ell\ge 
\gga(\beta,m)\ell\right\}.
\end{eqnarray}
Moreover for $0<m\le c_1$, we let $\tilde \bbP$ be the law obtained from $\bbP$
 shifting the distribution of 
$\go_1,\cdots,\go_\ell$  so that
\begin{eqnarray}
  \label{eq:to0}
\tilde \bbE[\omega_j]=8\frac{\gga(\beta,m)}{\beta m}.  
\end{eqnarray}
Note that 
\begin{equation}
  \label{eq:newmay}
  \lim_{m\to 0} {\gga(\beta,m)}/{ m}\,= \, 0,
\end{equation}
since otherwise,
 by convexity of $\gga(\beta,\cdot)$,
one has $\gga(\beta,m)\ge \epsilon m$ for for some $\epsilon>0$ and every $m$,
that is $ \tf(\beta,h_c(\gb)-\epsilon)=0$, which is in contrast 
with the definition of $h_c(\gb)$.

\smallskip
We now show that 
the event 
\begin{equation}
E:= \left\{ \go:\, \left(\go_1, \ldots, \go_\ell\right) \in  A_{\ell,m,\varepsilon}
\right\},
\end{equation}
becomes $\tilde \bbP$--typical
 for $\ell$ large. 

We first observe that,
thanks to the constraint $\cN_\ell/\ell\ge m-\varepsilon$, and assuming that $\varepsilon\le m/2$,
one has
\begin{eqnarray}
\label{uffa}
  \tilde\bbE\left(\frac 1 \ell{\log \hat Z_{\ell,\go}(\beta;m,\varepsilon)}-h_c(\beta)m
\right)\, \ge\,  4\gga(\beta,m)+
\bbE\left(\frac 1 \ell {\log \hat Z_{\ell,\go}(\beta;m,\varepsilon)}- h_c(\beta)m
\right).
\end{eqnarray}
By  \eqref{eq:def_phi} and
\eqref{eq:assurdo}, there
exist  $\varepsilon_0(m)>0$ and $\ell_0(\varepsilon,m)\in s\N$ such that
\begin{eqnarray}
  \tilde\bbE\left(\frac 1 \ell \log \hat Z_{\ell,\go}(\beta;m,\varepsilon)-h_c(\beta)m
\right)\ge 2\gga(\beta,m),
\end{eqnarray}
for $\varepsilon\le \varepsilon_0, \ell\ge \ell_0$. This in turn implies that 
\begin{eqnarray}
  \label{eq:PtildeE}
  \tilde\bbP\left(
  %\left(\go_1,\cdots,\go_\ell\right)\in A_{\ell,m,\varepsilon}
  E \right)\ge \tilde\bbP\left(
\log \hat Z_{\ell,\go}(\beta;m,\varepsilon)-\tilde\bbE \log \hat Z_{\ell,\go}(\beta;m,\varepsilon)\ge -
\gga(\beta,m) \ell
\right),
\end{eqnarray}
which is greater than, say, $1/2$ for $\ell$ sufficiently large, 
since $\gamma(\beta,m)>0$ and (recall \eqref{eq:def_phi} and discussion following that formula) 
\begin{eqnarray}
  \label{eq:convprob}
  \frac1 \ell {\log \hat Z_{\ell,\go}(\beta;m,\varepsilon)}-\frac1\ell\tilde\bbE\, {\log \hat Z_{\ell,\go}(\beta;m,\varepsilon)}
\stackrel{\ell\to \infty}{\longrightarrow}0,
\end{eqnarray}
in $\tilde\bbP$--probability. 
\medskip

The price of shifting $\bbP$ to $\tilde \bbP$ is directly estimated by using 
the assumption {\bf C2}, cf.
 \eqref{eq:exentropia}, and recalling \eqref{eq:newmay}: assuming that $m\le c_1(\beta)$, with 
$c_1(\beta)$  sufficiently small, one obtains the 
estimate
\begin{eqnarray}
  \label{eq:entropia}
  {\rm H(\tilde \bbP|\bbP):=\bbE \left( \frac{d\tilde\bbP}{d\bbP}\log \frac{d\tilde\bbP}{d\bbP}\right)}\le 
64R\left(\frac{ \gga(\beta,m)}{\beta m}\right)^2\ell,
\end{eqnarray}
and by applying the relative entropy inequality
\begin{eqnarray}
  \label{eq:B23}
  \log\left(\frac{\bbP(E)}{\tilde \bbP(E)}\right)
  \ge -\frac1{\tilde \bbP(E)}\left({\rm H}(\tilde \bbP|\bbP)+e^{-1}\right),
\end{eqnarray}
we obtain
\begin{equation}
  \label{eq:PA}
\tilde p:=  \bbP\left(
%\left(\go_1,\cdots,\go_\ell\right)\in A_{\ell,m,\varepsilon}
E \right)\, \ge\, 
\frac12\exp\left(-128R \left(\gga(\beta,m)/(\beta m)\right)^2
\ell\right),
\end{equation}
for large $\ell$. 

\smallskip
\begin{rem}\rm
\label{rem:Hineq}
Inequality \eqref{eq:B23}  holds whenever the measures $\bbP,\tilde\bbP$ are absolutely
continuous with respect to each other and for every event $E$ of nonzero 
measure. It  is a simple consequence of 
Jensen inequality: since $r \log r\ge -e^{-1}$ for every $r>0$, one has
\begin{equation}
  \label{eq:b23}
  \begin{split}
 \log\left(\frac{\bbP(E)}{\tilde \bbP(E)}\right)
 \,=&\, \log \tilde 
\bbE\left(\left.\frac{\dd\bbP}{\dd \tilde \bbP}\right| E\right)\ge \tilde 
\bbE \left(\log \left(\left.\frac{\dd\bbP}{\dd \tilde \bbP}\right)\right| E\right)
\\
\,=&\, -\frac1{\tilde\bbP(E)}\bbE \left(\frac{\dd\tilde\bbP}{\dd \bbP}\log \frac{\dd\tilde\bbP}{\dd \bbP} \ind_{E}\right)
\ge -\frac1{\tilde\bbP(E)}\left(
\bbE \left(\frac{\dd\tilde\bbP}{\dd \bbP}\log \frac{\dd\tilde\bbP}{\dd \bbP}\right)+e^{-1}\right),
\end{split}
\end{equation}
which is just \eqref{eq:B23}.
\end{rem}

\medskip

We now apply 
 an energy--entropy argument similar to that of \cite{cf:BG} which, in the present 
case, roughly consists in selecting only those polymer trajectories which visit the rare stretches where 
the disorder configuration is such to 
produce a sufficiently large positive fluctuation of the partition function. 
Of course the precise definition of these rare stretches is directly
related to the event $E$. 
This selection strategy 
gives a lower bound on the free
energy, which implies \eqref{eq:stimaphi}.
More precisely, we consider a system of length $k\ell$, with $k\in \N$, $\ell\in s\N$ 
and we divide it into blocks 
$B_j=\{j \ell+1, j\ell+2,\cdots,(j+1)\ell\}$ of length $\ell$, with $j=0,\cdots,k-1$. For a given realization of 
$\go$, we denote by $\mathcal I(\go)$ the ordered set of nonnegative integers
\begin{eqnarray}
  \label{eq:I}
  \mathcal I(\go) =\{j_1,\cdots,j_{|\mathcal I(\go)|}\}:=
\{0\le j\le k-1: \left(\go_{j \ell+1},\cdots,\go_{(j+1)\ell}\right)\in A_{\ell,m,\varepsilon}\},
\end{eqnarray}
where it is understood that $\varepsilon\le \varepsilon_0(m)$, $\ell\ge \ell_0(m,\varepsilon)$ and $m\le c_1$
as above, so that
\eqref{eq:PA} holds.
We bound the partition function $Z_{N,\go}$ below by
inserting in the average over the  paths the constraint that $S_i\ne0$ whenever $i\in B_j$ with 
$j\notin\mathcal I(\go)$,
that $S_i=0$ whenever $i=j\ell$ or $i=(j+1)\ell$ with $j\in \mathcal I(\go)$, and that, for every
$j\in \mathcal I(\go)$,
\begin{eqnarray}
\frac{|\{i\in B_j: S_i=0\}|}\ell\in [m-\varepsilon,m+\varepsilon].
\end{eqnarray}
\begin{figure}[h]
\begin{center}
\leavevmode
\epsfysize =2 cm
\psfragscanon
\psfrag{0}[c][l]{$0$}
\psfrag{n}[c][l]{ $n$}
\psfrag{b1}[c][]{ $B_{1}$}
\psfrag{b2}[c][]{ $B_{2}$}
\psfrag{b3}[c][]{ $B_{3}$}
\psfrag{b4}[c][]{ $B_{4}$}
\psfrag{b5}[c][]{ $B_{10}$}
\psfrag{b6}[c][c]{ $B_{17}$}
\psfrag{dots}[c][c]{ $\dots$}
\psfrag{L0}[c][c]{ $L_0$}
\psfrag{L1}[c][c]{ $L_1$}
\psfrag{L2}[c][c]{ $L_2$}
\psfrag{L3}[c][c]{ $L_3$}
\psfrag{N}[c][l]{ $k\ell$}
\psfrag{Sn}[c][l]{$S_n$}
\epsfbox{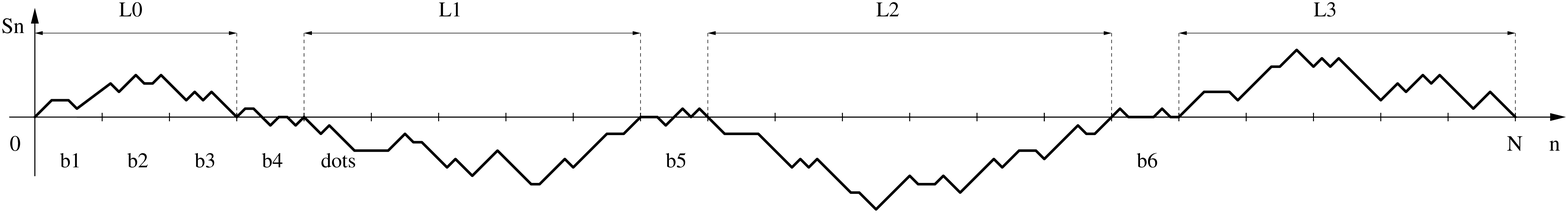}
\end{center}
\caption{\label{fig:figura} A typical trajectory contributing to the lower bound \eqref{eq:lbZ}.
In this example $k=22$, $\ell=8$, $\mathcal I(\go)=\{4,10,17\}$. Note that $S_n\ne 0$ for $n$ in  a 
block $B_j$ with $j \not\in \{4,10,17\}$, except at the boundary with a block $B_j$ with 
$j \in \{4,10,17\}$,
 since by construction $S_n=0$ at the boundaries of these blocks.
 Inside $B_{j}$,  $j \in \mathcal I(\go)$ ,  the path moves with relative freedom, but it is
 bound to touch
the line $S=0$ approximately $m\ell$ times (see the text for the choice of $m$).}
\end{figure}
In this way, recalling the definition of $A_{\ell,m,\varepsilon}$ and assuming  that
$\varepsilon\le \gga(\beta,m)/(2|h_c(\beta)|)$, one obtains 
\begin{eqnarray}
  \label{eq:lbZ}
  \frac1{k\ell}{\log Z_{k\ell,\go}(\beta,h_c(\beta))}\ge \frac 1k|\mathcal I(\go)|
\frac{\gga(\beta,m)}2+\frac1{k\ell}
\log \prodtwo{r=0,\cdots,|\mathcal I(\go)|:}{L_r\ne0} K(L_{r}),
\end{eqnarray}
where we recall that $K(n)$ is the $\bP$--probability that first return to $0$ of an excursion of the free
process occurs at step $n$, as in
 \eqref{eq:ret_times}, while  $L_r$'s, the (possibly vanishing) lengths of the 
excursions of the process between two blocks $B_i,B_{i'}$ with $i,i'\in \mathcal I(\go)$, 
are defined as
\begin{eqnarray}
  \label{eq:Lj}
  L_r=\ell|j_{r+1}-j_r-1|,
\end{eqnarray}
with the convention that $j_0=-1$, $j_{|\mathcal I(\go)|+1}=k$,
%and that we set $\log K(L_r)\equiv 0$ whenever $L_r=0$  
see Figure  \ref{fig:figura}.
Taking the expectation with respect to the disorder and using  \eqref{eq:RW}, one obtains then
\begin{eqnarray}
  \frac1{k\ell}\bbE \log Z_{k\ell,\go}(\beta,h_c(\beta))\, 
  \ge\,
  \frac{\gga(\beta,m)}2   \tilde p-
\frac{2\alpha}{k\ell}\bbE \sum_{r=0}^{|\mathcal I(\go)|}\ind_{L_r\ne0} \log L_r,
\end{eqnarray}
for $\ell$ sufficiently large.
At this point, as in \cite{cf:BG}, one uses Jensen's inequality and the concavity of the logarithm to
get
\begin{equation}
  \label{eq:jensen}
  \begin{split}
   \frac1{k\ell}\bbE \log Z_{k\ell,\go}(\beta,h_c(\beta))
   \,&\ge \,\frac{\gga(\beta,m)}2  \tilde p 
-\frac{2\alpha}{k\ell}\bbE\left[\left(|\mathcal I(\go)|+1\right)\log \left(\frac{\sum_{r=0}^{|\mathcal I(\go)|}
L_r}{|\mathcal I(\go)|+1}\right)\right]\\
&\ge\, \frac{\gga(\beta,m)}2  \tilde p
-2\alpha\bbE\left[\frac{\left(|\mathcal I(\go)|+1\right)}{k\ell}\log \left(\frac{k\ell}{|\mathcal I(\go)|+1}\right)\right].
\end{split}
\end{equation}
Since the disorder variables in the distinct blocks $\left\{B_j\right\}_j$ are independent,  the 
law of large numbers implies
\begin{eqnarray}
  \frac{|\mathcal I(\go)|}{k}\stackrel{k\to\infty}{\longrightarrow}{\tilde p}
\end{eqnarray}
$\bbP(\dd \go)$--a.s.,
 so that, recalling \eqref{eq:PA},
\begin{equation}
\label{finale}
\begin{split}
0=  \tf(\beta,h_c(\beta))\, &=\, \lim_{k\to\infty}\frac1{k\ell}\bbE \log Z_{k\ell,\go}(\beta,h_c(\beta))
\\ & \ge \,\tilde p\frac {\gga(\beta,m)}2-\frac{2\alpha\tilde p}\ell\log\frac\ell{\tilde p}\,
\ge\, \tilde p\left(\frac {\gga(\beta,m)}2-256\alpha R\left(\frac{\gga(\beta,m)}{\beta m}\right)^2
-2\alpha\frac{\log (2\ell)}\ell\right).
\end{split}
\end{equation}
Since $\ell$ is arbitrary, one obtains
\begin{eqnarray}
\label{eq:limiti}
\gga(\beta,m)\ge \frac{\beta^2 m^2}{512 \alpha R}
\end{eqnarray}
which is the desired inequality and \eqref{eq:stimaphi}, provided that $c_2$ in \eqref{eq:assurdo} 
satisfies $c_2<(512 R)^{-1}$.

\hfill $\stackrel{\text{Theorem  \ref{th:stimaphi}, {\bf C2}}}{\Box}$

\begin{rem}\rm 
It is easy to check that, in the Gaussian case $\go_1\sim \mathcal N(0,1)$, Theorem \ref{th:smooth_cont}
holds with $c(\beta)=c_3 \beta^{-2}$, $c_3$ a suitable positive constant. Indeed, in this case
 the estimate \eqref{eq:exentropia} holds for every $x\in \R$ (and $R=1/2$) and therefore one can take $c_1=1$ in 
the proof of Lemma \ref{th:stimaphi}, so that \eqref{eq:limiti} implies \eqref{eq:stimaphi} 
with $\cone(\beta)=c_4\beta^2$. Of course, the same is true, up to constants, for every $\bbP$ such that  inequality \eqref{eq:exentropia}
holds uniformly for $x\in \R$.
\end{rem}
\medskip

\noindent
{\it Proof of Theorem \ref{th:stimaphi} under assumption {\bf C1}}.
The proof proceeds as in case {\bf C2}, up to  the definition of the law $\tilde \bbP$,
but in this case
the law obtained by shifting $\bbP$  in general
has an infinite entropy with respect to $\bbP$. Therefore, in this case we define $\tilde \bbP$ rather by tilting 
the law of the first $\ell$ variables:
\begin{eqnarray}
  \label{eq:tildebbP}
 \frac{ \dd \tilde \bbP}{\dd \bbP}( \go)\, =\,  
 \frac1{z^{\ell}}
 \exp\left({u \sum_{n=1}^\ell \omega_n}\right),
\end{eqnarray}
where $u\ge 0$ will be  chosen later and $z=z(u)=\bbE\, e^{u\go_1}$.
Let, for $\ell\in s\N$,
\begin{eqnarray}
  \label{eq:psi}
  \psi_\ell(u)=\frac1\ell\tilde \bbE \,
\log \hat Z_{\ell,\go}(\beta;m,\varepsilon)
\end{eqnarray}
and observe that 
$$
\psi_\ell(0)=\frac1\ell\bbE \,\log \hat Z_{\ell,\go}(\beta;m,\varepsilon).
$$
Then, one has 
\medskip

\begin{lemma}
\label{lemma_psi}
There exist $u_0(\beta)>0$ and $c_0(\beta)>0$, possibly depending on $\bbP$, such that
for every $0<\beta<\infty$, $1\le \alpha<\infty$  the following holds:  for every $m\in (0,1]$, if $\varepsilon\le m/2$ 
and $0\le u\le u_0$ we have
\begin{eqnarray}
  \label{eq:stimapsi}
  \psi_\ell(u)-\psi_\ell(0)\ge  c_0\beta m u.
\end{eqnarray}
\end{lemma}
\medskip

Lemma \ref{lemma_psi} will be proven below. 
To proceed with the proof of Theorem \ref{th:stimaphi}, we choose $u= 4\gamma(\beta,m)/(\beta m c_0)$ 
and notice that $u$ is
 certainly smaller than $u_0$ if $m\le c_1(\gb)$ with $c_1$ sufficiently small (see 
 \eqref{eq:newmay}). Then, choosing $0< m\le c_1$ and  $\varepsilon\le m/2$,  \eqref{eq:stimapsi} implies that
\eqref{uffa} is valid also in the present case. On the other hand, it is immediate to verify that 
 \eqref{eq:entropia} still holds, with 
$R$ replaced by some $R'(\beta,M)$ and $\gamma(\beta,m)$ by $\gamma(\gb,m)/c_0$. 
The rest of the proof proceeds exactly as in the case {\bf C2}.

\hfill $\stackrel{\text{ Theorem   \ref{th:stimaphi}, {\bf C1}}}{\Box}$

\medskip

\noindent
{\it Proof of Lemma \ref{lemma_psi}}.
Note that
\begin{eqnarray}
  \label{eq:der_u}
  \partial_u  \psi_\ell(u)=\frac1{ \ell}\sum_{i=1}^\ell 
\tilde \bbE \left[\,\go_i\,
\log \hat Z_{\ell ,\go}(\beta;m,\varepsilon)\right]
-\ell \xi\,\psi_\ell(u),
\end{eqnarray}
where 
$\xi=\xi(u)=\tilde\bbE \,\go_1.$
The first term in the right--hand side of \eqref{eq:der_u} can be rewritten 
(with some abuse of notation) as 
\begin{multline}
\label{elimin}
\frac{\xi}{\ell}\sum_{i=1}^\ell
\tilde \bbE \, \left[\left.
\log \hat Z_{\ell,\go}(\beta;m,\varepsilon)\right|_{\go_i=0}\right]+\frac {\beta}{\ell}\sum_{i=1}^\ell 
\tilde \bbE \left[\,\go_i\, \int_{0}^{\go_i}\dd y
\left.\hat\bE_{\ell,\go}(\ind_{S_i=0}) \right|_{\go_i=y}\right]
\,
\\
=\ell \xi\psi_\ell(u)+
\frac {\beta}{\ell}\sum_{i=1}^\ell \tilde \bbE\,\left[ (\go_i-\xi) 
\int_0^{\go_i}\dd y\left.\hat\bE_{\ell,\go}(\ind_{S_i=0}) \right|_{\go_i=y}\right],
\end{multline}
where
\begin{eqnarray}
\hat\bE_{\ell,\go}(\cdot)=\frac{
\bE 
\left[ \cdot\,
\exp\left( \sum_{i=1}^{\cN_\ell}
\gb \go_{\tau_i}
\right)
\ind_{S_\ell=0}\ind_{(\cN_\ell/\ell)\in [m-\varepsilon,m+\varepsilon]}
\right]
}{\hat Z_{\ell ,\go}(\beta;m,\varepsilon)}
\end{eqnarray}
and obviously the first term in \eqref{elimin} cancels the second one in 
the right--hand side of \eqref{eq:der_u}.
Next, observe that the following identity holds:
\begin{eqnarray}
\label{eq:mgc}
 \left.\hat\bE_{\ell,\go}(\ind_{S_i=0}) \right|_{\go_i=y}=\frac{e^{\gb (y-y')}\left.\hat\bE_{\ell,\go}(\ind_{S_i=0})\right|_{\go_i=y'}}
{1+(e^{\gb (y-y')}-1)\left.\hat\bE_{\ell,\go}(\ind_{S_i=0})\right|_{\go_i=y'}}.
\end{eqnarray}
Now recall \eqref{eq:omegabound}, so that it is sufficient to consider
 $y$ and $y'$ such that $\vert y-y'\vert \le  2M$, and from \eqref{eq:mgc}
 we obtain
\begin{eqnarray}
e^{-2\gb M}\left.\hat\bE_{\ell,\go}(\ind_{S_i=0})\right|_{\go_i=y'} 
\le   \left.\hat\bE_{\ell,\go}(\ind_{S_i=0}) \right|_{\go_i=y}\le 
e^{+2\gb M}\left.\hat\bE_{\ell,\go}(\ind_{S_i=0})\right|_{\go_i=y'}.
\end{eqnarray}
We can use this inequality to bound below the last term in \eqref{elimin}. We have in fact
\begin{eqnarray}
\nonumber
  \label{eq:primoLB}
 \tilde \bbE\left[\, \go_i 
\int_0^{\go_i}\dd y\left.\hat\bE_{\ell,\go}(\ind_{S_i=0}) \right|_{\go_i=y}\right]\ge
e^{-2\gb M}\eta\,\tilde \bbE
\left.\hat\bE_{\ell,\go}(\ind_{S_i=0}) \right|_{\go_i=0}
\ge e^{-4\gb M}\eta\,\tilde \bbE\,
\hat\bE_{\ell,\go}(\ind_{S_i=0}),
\end{eqnarray}
where $\eta=\eta(u)=\tilde \bbE \,\omega_1^2$, while
\begin{eqnarray}
  \label{eq:secondoLB}
  \tilde \bbE\, 
\int_0^{\go_i}\dd y\left.\hat\bE_{\ell,\go}(\ind_{S_i=0}) \right|_{\go_i=y}\le Me^{2\gb M}\tilde \bbE\,
\hat\bE_{\ell,\go}(\ind_{S_i=0}) .
\end{eqnarray}
Therefore, recalling the constraint 
$(\cN_\ell/\ell)\in [m-\varepsilon,m+\varepsilon]$  in the definition of $\hat\bP_{\ell,\go}(\cdot)$,
one has the following lower bound:
\begin{eqnarray}
\label{LB}
 \partial_u  \psi_\ell(u)\ge \beta( m-\varepsilon)\eta
e^{-4M\gb}-\beta(m+\varepsilon)\xi M e^{2\gb M}.
\end{eqnarray}
Now choose $\varepsilon\le m/2$ and notice that
$\xi=u+O(u^2)$ for 
$u\ll 1$, while $\eta=1+O(u)$. 
Therefore, there exists $u_0(\beta,M)>0$ such that, for  $0\le u\le u_0$ 
and for every $m$, the following holds: 
\begin{eqnarray}
 \partial_u  \psi_\ell(u)\ge \beta m\frac{e^{-4M\gb}}{4}-2\beta m u Me^{2\gb M}
\ge \beta m \frac{e^{-4M\gb}}{8}=: c_0 \beta m.
\end{eqnarray}
An integration in $u$ concludes the proof of \eqref{eq:stimapsi}.

\hfill $\stackrel{\text{ Lemma \ref{lemma_psi}}}{\Box}$

\subsection{The copolymer case}
\label{sec:prova_cop}
In order to prove Theorem \ref{th:smooth_cont}  for the copolymer model
\eqref{eq:Hcop}, it is convenient to start from the observation that one can rewrite the limit free energy as
\begin{eqnarray}
  \tf (\beta,h)=\lim_{N\to\infty}\frac1N \bbE \log \bE\left[ \exp\left(-\sum_{n=1}^N(\gb \go_n+ h)\Delta_n\right)
\ind_{S_N=0}\right],
\end{eqnarray}
where $\Delta_n=\ind_{\sign(S_n)=-1}$. Comparing this expression for the 
copolymer free energy with  \eqref{eq:RNS}, it is clear that in the present case
the role of $\ind_{S_n=0}$  is played by  $\Delta_n$.  The proof then proceeds exactly like in the pinning
case, with the only differences that in the definition \eqref{eq:def_phi} of $\phi(\beta,m)$
the constraint $(\cN_N/N)\in [m-\varepsilon,m+\varepsilon]$ has to be replaced by
\begin{eqnarray}
  \frac{|\{1\le n\le N: \Delta_n=1\}|}N\in [m-\varepsilon,m+\varepsilon]
\end{eqnarray}
and that, in the energy--entropy argument, 
the path is required to satisfy $S_i>0$ (and not just $S_i\ne0$)
whenever $i\in B_j$ with $j\notin \mathcal I(\go)$.
This implies that $K(L)$ in  \eqref{eq:lbZ} has to be replaced by $K(L)/2$, which has the effect of adding a negative term of 
order $O(\tilde p/\ell)$ in the lower bound \eqref{finale}, which is negligible for $\ell$ sufficiently large.

\hfill $\stackrel{\text{ Theorem   \ref{th:stimaphi}, copolymer}}{\Box}$

\appendix 

\section{The non--disordered pinning model}
\label{app:b0}
For $\gb=0$,  the free energy \eqref{eq:free_energy_p}
may be identified explicitly with the following procedure. 
First we consider the equation
\begin{equation}
\label{eq:b}
\sum_{n\in \N} K(n)\exp(-bn)= \exp(h), 
\end{equation}
and we look for a solution $b>0$, which exists only if $ \sum_{n\in\N} K(n)>\exp(h)$,
in which case
it is unique. 
Then  if we set $\tilde K (n) := \exp(-h-bn)K(n)$, $\tilde K (\cdot)$ is a
discrete probability density and  
 one can write
\begin{equation}
\begin{split}
\bE\left[ \exp\left( -h \cN_N \right); \, \tau_{\cN_N}=N \right]\, &=\,
\sum_{j=1}^N \sumtwo{(\ell_1, \ldots , \ell_j )\in \N^j:}{
\sum_{i=1}^j \ell_i =N} \prod_{i=1}^j \exp(-h) K(\ell_i)
\\
&= \, \exp\left( b N\right) 
\sum_{j=1}^N \sumtwo{(\ell_1, \ldots , \ell_j )\in \N^j:}{
\sum_{i=1}^j \ell_i =N} 
\prod_{i=1}^j \tilde K(\ell_i) =:  \exp\left( b N\right) G_N,
\end{split}
\end{equation}
and one easily sees that $G_N$  is the probability that the random walk 
which starts at $0$ and takes positive integer IID jumps with law $\tilde K(\cdot)$
hits the site $N$. It is a classical fundamental result of renewal theory that $\lim_N G_N = 1/
\sum_{n\in \N} n \tilde K(n)$
\cite[Ch. XIII]{cf:Feller1}. This of course implies that $b= \tf(0,h)$.
On the other hand, if \eqref{eq:b} admits no positive solution, by proceeding 
as in \eqref{eq:b} and by setting simply $\tilde K (n) := \exp(-h)K(n)$, so that
$\tilde K( \cdot)$ is a sub--probability density, one easily sees that $\tf (0, h)=0$.
So equation equation \eqref{eq:b} contains all the information about the free
energy. 

Let us then observe  that
 \eqref{eq:b} has a positive solution if and only if 
$h<\log (1-K(\infty))$ and therefore $h_c(0)=\log (1-K(\infty))\le 0$ ($h_c(\gb)$ being defined before
Remark \ref{rem:free}).
The behavior at criticality can be extracted from
\eqref{eq:b} in a rather straightforward way too, but of course we need to make precise
the requirement on $K(\cdot)$ beyond the lower bound
\eqref{eq:RW}:
\begin{itemize}
\item {\sl The case $\gS:=\sum_{n\in \N} n K(n)<\infty$}. This is a necessary and sufficient
condition for the transition to be of first order.
More precisely, for $\gS \in (0, \infty]$
\begin{equation}
\label{eq:1ord}
\tf \left( 0, h\right)\, =\, \frac{\exp(h_c(0))}{\gS} \left(h_c(0) -h\right) +o\left(
\left(h_c(0) -h\right)\right), \ \ \text{ for } h \nearrow h_c(0).
\end{equation}
Formula \eqref{eq:1ord} follows since by Dominated Convergence if $\gS< \infty$
then 
\begin{equation}\sum_{n\in \N } K(n)\exp(-bn) =\exp(h_c(0))- b \gS +o(b)
\end{equation}
 for $b \searrow 0$,
while by a direct estimate $\lim_{b \searrow 0}b^{-1}\sum_{n\in \N} K(n)[1-\exp(-bn)]=\infty$
if $\gS=\infty$.
Note that the condition $\gS<\infty$ holds, in particular, in the case
where 
\begin{equation}
\label{eq:e.g.}
K(n) \, = \, L(n)/n^\ga   
\end{equation}
for large $n\in s\N$, with $\ga>2$ and $L(\cdot)$ a function varying slowly 
at infinity, i.e., a positive function such that $\lim_{r\to \infty} L(xr)/L(r)=1$ for every 
$x>0$ (see \cite{cf:regular} for more details on slowly varying functions).

\item {\sl The case $\gS=\infty$}. We
set $\overline{K} (n)= \sum_{j=n+1}^\infty K(n)$ (by this we mean the sum
over $j \in \N$ such that $j>n$: $K(\infty)$ is not included in the sum)
and assume that 
$$\sum_{j=1}^N \overline{K} (j) = \widehat L(N) N^{2-\ga}, $$
for some function $\widehat L(\cdot)$ which is slowly varying at infinity. This is true, in particular, in the case \eqref{eq:e.g.}.
By the easy (Abelian) part of the classical Tauberian Theorem
\cite[Ch.XIII.5, Theorem 2]{cf:Feller2}
we have that $\sum_{n=0}^\infty \exp(-bn)\overline{K}(n)=c^\prime b^{\ga -2}\widehat L(1/b) (1+o(1))$ as $b \searrow 0$,
with $c^\prime= c^\prime(\ga)>0$.
Therefore
\begin{equation}
\label{eq:ipTaub}
\begin{split}
\sum_{n=1}^\infty e^{-bn}K(n)\, &=\, \overline K (0) +
 \left( e^{-b}-1 \right) \sum_{n=0}^\infty e^{-bn}\overline{K}(n)\\
  &=\, \overline K (0) - c^\prime b^{\ga -1}\widehat L(1/b) (1+o(1)),
  \end{split}
\end{equation} 
and  
\begin{equation}
\label{eq:holder}
\tf (0,h)\,= \left(h_c(0) -h\right)^{1/(\ga -1)} \tilde L\left(1/(h_c(0) -h)\right), \ \ 
\text{ for } h \nearrow  h_c(0),
\end{equation}
with $\tilde L(\cdot)$
a slowly varying function (see \cite[(1.5.1) and Theorem 1.5.12]{cf:regular}).
It is therefore clear that the transition is of second order for $\ga \in (3/2,2]$
(we emphasize, for the case $\alpha=2$, that we are assuming that $\gS=+\infty$)
and it is of higher order for $\ga<3/2$.
The value $\ga=3/2$ is borderline and the order of the transition depends then on 
the slowly varying function $\widehat{L}(\cdot)$. In the case of one dimensional symmetric random walks
with IID increments taking values in $\{-1,0,+1\}$,  $\ga=3/2$ and $L(n)=c (1+o(1))$ for $n$ large,
$c$ a positive constant,    
 and  therefore the transition  is really of second order and not higher.

\end{itemize}

\medskip

\section*{acknowledgments} We would like to thank Bernard Derrida, Thomas Garel, C\'ecile Monthus and David Mukamel
for interesting discussions.
This research has been conducted in the framework of the GIP--ANR project 
 JC05\_42461
({\sl POLINTBIO}).


\begin{thebibliography}{15}  


 
  

\bibitem{AizWeh90}
M. Aizenman and J. Wehr, \textit{
Rounding effects of quenched randomness on first--order phase transitions},
Comm. Math. Phys. {\bf 130} (1990),  489--528. 

\bibitem{cf:AS}
 K. S. Alexander and V. Sidoravicius,
 \textit{Pinning of polymers and interfaces by random potentials},
 preprint (2005). Available on: arXiv.org e-Print archive: math.PR/0501028 


\bibitem{cf:regular} N. H. Bingham, C. M. Goldie and J. L. Teugels, \textit{Regular Variation},  Cambridge University Press, Cambridge, 1987.

\bibitem{cf:BG}
T. Bodineau and G. Giacomin, \textit{On the localization transition   of random copolymers near selective interfaces},
J. Statist. Phys. {\bf 117} (2004), 801--818.

\bibitem{cf:BdH} E. Bolthausen and F. den Hollander,    
 \textit{Localization transition for a polymer near an interface},   
 Ann. Probab.  {\bf 25}  (1997),  1334--1366.   
   
 \bibitem{cf:BK}  
A.   Bovier and C.  K\"ulske, 
\textit{There are no nice interfaces in $(2+1)$--dimensional SOS models in random media}, 
J. Statist. Phys. {\bf 83} (1996),  751--759.
 
\bibitem{cf:CCFS} J.~T.~Chayes, L.~Chayes, D.~S.~Fisher and T.~Spencer, {\em
Correlation Length Bounds for Disordered Ising Ferromagnets}, Commun. Math. Phys. {\bf 120}, 501--523
(1989).

 \bibitem{cf:C}
 B. Coluzzi, \textit{Numerical study on a disordered model for DNA denaturation transition},
 arXiv:cond--mat/0504080v1 (2005).
 
 \bibitem{cf:CH}
 D. Cule and T. Hwa, \textit{Denaturation of heterogeneous DNA}, Phys. Rev. Lett. {\bf 79} (1997), 2375--2378.
   
   
\bibitem{cf:Derrida} B. Derrida, V. Hakim and J. Vannimenius, \textit{Effect of disorder on 
two--dimensional wetting}, J. Statist. Phys. {\bf 66} (1992), 1189--1213. 

\bibitem{cf:Feller1}  
W.~Feller, \textit{An introduction to probability theory and its applications}, Vol. I,  
 Third edition, John Wiley \& Sons, Inc.,   
New York--London--Sydney, 1968.  

\bibitem{cf:Feller2}  
W.~Feller, \textit{An introduction to probability theory and its applications}, Vol. II,  
 Second edition, John Wiley \& Sons, Inc.,   
New York--London--Sydney, 1971.  


\bibitem{cf:Fogacs} G. Forgacs, J.~M.~Luck, Th.~M.~Nieuwenhuizen and H.~Orland, 
\textit{Wetting of a Disordered Substrate: Exact Critical behavior  in Two Dimensions}, 
Phys. Rev. Lett. {\bf 57} (1986), 2184--2187.

\bibitem{cf:GHLO}
T. Garel, D. A. Huse, S. Leibler and H. Orland,
\textit{Localization transition of random chains at interfaces},
Europhys. Lett. {\bf 8} (1989), 9--13.

\bibitem{cf:GM2}
T. Garel and C. Monthus,
\textit{Numerical study of the disordered Poland--Scheraga model of DNA denaturation},
 J. Stat. Mech., Theory and Experiments (2005), P06004.

\bibitem{cf:G}  
G.~Giacomin, \textit{Localization phenomena in random polymer models}, preprint (2004). 
\\
Available
online: http://www.proba.jussieu.fr/pageperso/giacomin/pub/publicat.html

\bibitem{cf:GT} G.~Giacomin and  F.~L.~Toninelli, \textit{Estimates on path delocalization for copolymers 
at selective interfaces}, Probab. Theor. Rel. Fields {\bf 133} (2005), 464--482.

\bibitem{cf:GTloc} G.~Giacomin and F.~L.~Toninelli, \textit{The localized phase of disordered copolymers with 
adsorption}, preprint (2005). Available on: arXiv.org e--Print archive: math.PR/0510047

\bibitem{cf:Harris} A.~B.~Harris, \textit{Effect of random defects on the critical behaviour of Ising models}, 
J. Phys. C {\bf 7} (1974), 1671--1692.

\bibitem{cf:Imry-Ma} Y. Imry and S.--K. Ma, \textit{Random--Field Instability of the Ordered State of Continuous Symmetry},
Phys. Rev. Lett. {\bf 35} (1975), 1399--1401.

\bibitem{cf:KMP}
Y. Kafri, D. Mukamel and L. Peliti,
\textit{Why is the DNA denaturation transition first order?}, Phys. Rev. Lett.
{\bf 85} (2000), 4988--4991.

\bibitem{cf:King} J. F. C. Kingman, \textit{Subadditive ergodic theory}, Ann. Probab. {\bf 1} (1973),
882--909.

\bibitem{cf:Monthus} C. Monthus, \textit{On the localization of   
random heteropolymers at the interface between two selective solvents},  
Eur. Phys. J. B {\bf 13} (2000), 111--130.   
  
\bibitem{cf:Nicolas} N. Petrelis, \textit{Polymer pinning at an interface}, 
preprint (2005). Available on: arXiv.org e--Print archive:  math.PR/0504464
  
\bibitem{cf:Sinai}
Ya. G. Sinai, \textit{A random walk with a random potential},  Theory Probab. Appl. {\bf 38} (1993),  382--385.

\bibitem{cf:SW} C. E. Soteros and S. G. Whittington,  \textit{The statistical mechanics of random copolymers}, J. Phys. A: Math. Gen. {\bf 37} (2004), R279--R325.

\bibitem{cf:TangChate} L.--H. Tang and H. Chat\'e, \textit{Rare--Event Induced Binding Transition of 
Heteropolymers}, Phys. Rev. Lett. {\bf 86} (2001), 830--833.

\bibitem{cf:TM}
  A. Trovato  and A. Maritan,
\textit{A variational approach to the localization transition of heteropolymers at interfaces},
Europhys. Lett. {\bf 46} (1999), 301--306.

\end{thebibliography}
\end{document}